**PAPER • OPEN ACCESS**

# Assessment various control methods a digital copy of enterprise by integral indicator

View the article online for updates and enhancements.

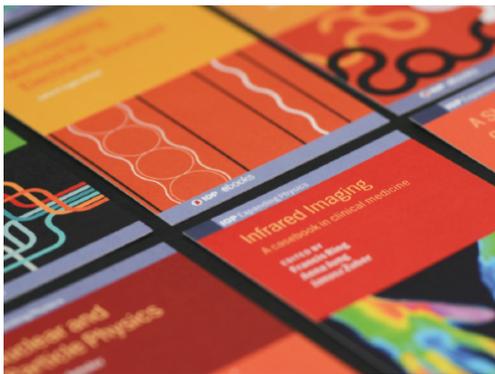

# Assessment various control methods a digital copy of enterprise by integral indicator

**S N Masaev**[1,2]

[1] Siberian Federal University, pr. Svobodnyj, 79, Krasnoyarsk, 660041, Russia
[2] Control Systems LLC, 86 Pavlova Street, Krasnoyarsk, 660122, Russia

E-mail: faberi@list.ru

**Abstract**. The difficulty of assessing the state lies in a little predictable change in the dimension of a dynamic system under the influence of internal changes and environmental parameters. In the work, the state of such a system is estimated by the method of integral indicators. The application of the method of integral indicators allowed us to evaluate the activity of an enterprise. In the present work, the method of integrated indicators is used to assess the control of a digital copy (enterprise). Using the author's complex of programs, the activity of an enterprise (digital copy) is modeled as a dynamic non-stationary, in structure, system of 1.2 million values with spaces identified at each time step: output data, control data, and environmental parameters. The research showed significant changes in the values of the integral indicator characterizing the state of a dynamic system when implementing different regimes of control for a woodworking enterprise and a construction company. In the first and second examples, an underestimation of the invisible material flows of production processes of "shadow" businesses-processes were found to be 4.5 billion rubles per year. In the third example, the loss of the enterprise from stop businesses-processes for 4 years will amount to 188.8 million rubles.

## 1. Introduction
A large number of approaches to assess the regimes of control enterprise have been identified [1-4]. However, no universal method for a unified assessment of the control modes of a digital copy of an enterprise has been found [5] Purpose: to estimate various control methods of digital copy by integral indicator [6-14].

## 2. Method
The description of a system, parameters and control loop is described in detail in the previous work [11-13]. Dynamic and non-stationary, in structure, the system is given by the equation

$$y(t)=A(t)x(t)+B(t)u(t)+v(t), \qquad (1)$$

$$y(t)=G(x(t),u(t),v(t))$$

$$x(t) \in R^n,\ u(t) \in R^m,\ v(t) \in R^l,\ y(t) \in R^k$$

where $x(t) \in R^n$ is the vector of state signals of the object of research (income/expense items), $u(t) \in R^m$ is the vector of control actions at the level systems (planning of income/expense items), $v(t) \in R^l$ is the vector of environmental perturbations, $y(t) \in R^k$ is the observation vector, and $y(t)=Hx(t+1)$, $A=[a_{ij}]$ –

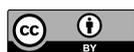

$N \times N$ is the matrix that determines the development speed and structure of the interaction of the parameters of the $x(t)$ system, $B=[b_{ij}] - N \times M$ is the matrix determining the structure of the control action of $u(t)$ on the system. The analysis of the system at the time $t$ is performed on $x(t)$ for $k$ of the previous measures. The parameter $k$ is the length of the time series segment. Then we have the matrix $X_k(t)$. The value of $k$ can be used to identify the influence of external factors on the system [15]. Following the method of integral indicators. We get the correlation matrix $R_k(t)$ with the correlation coefficients $r_{ij}(t)$:

$$R_k(t) = \frac{1}{k-1} X_k^{oT}(t) X_k^o(t) \qquad (2)$$

$$r_{ij}(t) = \frac{1}{k-1} \sum_{l=1}^{k} x_k^{oi}(t-l) x_k^{oj}(t-l), \quad i,j = 1,\ldots,n, \qquad (3)$$

and calculate the integral indicator

$$G = \sum_{t=1}^{T=\max} \sum_{i=1}^{n} G_i(t), \qquad (4)$$

$$G_i(t) = \sum_{j=1}^{n} |r_{ij}(t)|.$$

Goal control $J=J(W) \to \inf_{w \in \Omega^*}$, where $t$ the procedure for predicting the state of a dynamic system in different clock cycles, $W=[w_{ij}]=T \times N$ is the matrix that determines the structure of the chosen control method, in the general form $J$ - is the goal parameter for the states of the system for various t, which determines the admissible set of matrices of the selected control method $\Omega^*$, taking into account all the requirements for the process dynamics and the language syntax of the selected control method $u=W(t)y$.

The implementation of the method is performed in the author's complex of programs.

## 3. Calculation algorithm

Author's software package consists of four separate programs: 1. A software package for calculating the economic model of the functioning of an enterprise engaged in the procurement and deep processing of various types of wood under certain scenarios of market development and strategy (RosPatent Certificate on registration of a computer program No. 2013614410). The program not only simulates the activities of an enterprise in the forest industry, it is also used to model any economic object (enterprise) as a system, to determine the structure of the interaction of variables, to simulate a managerial decision, the structure of the interaction of variables in a control action, and to simulate the frequencies of the external environment acting on the system. 2. Assessment of the achievement of the chosen strategy of the enterprise by universal indicators of business functions of production processes and control processes (RosPatent Certificate on registration of a computer program No. 2017616973). Used to link the system parameters with the selected control method. Allows you to choose a convenient form of data interpretation according to the syntax of the control method. 3. The software package for evaluating the effectiveness of a managerial decision (RosPatent Certificate on registration of a computer program No. 2008610295). The method of R. Bellman is used to assess the optimality of managerial decisions. 4. Automated calculation and filling out forms of economic evaluation of investment projects in accordance with the decision of the Council of the Administration of the Krasnoyarsk Territory "On state support of investment activities 91-P" (RosPatent Certificate on registration of a computer program No. 2017616970). It is used to model the optimal regime of tax rates by the subject of the Russian Federation (region, region, etc.) to stimulate the development of enterprises of residents of the SEZ (Special Economic Zone).

This software package is used to perform the experiment algorithm. The step number in the algorithm coincides with the number of the program being executed from the author's program complex described above.

1 step. We load data on the actual activity of the economic object $x_i(t)$, control data for past periods of time and data characterizing the external environment. We identify the enterprise as a multidimensional dynamic system. We introduce a strategy at the enterprise. We calculate the integral indicator for all states of the system. We set the objective function. If the data characterize the object, then go to step 2, otherwise repeat 1 step.

2 step. In the loaded model, we select the form of representation of the control method. For example, PMBOK [16], budgeting, Cobb-Douglas formula, V-shaped project control model, quality management system (TQM), life cycle assessment method (products, enterprises, strategies), etc. Additional examples can be found in separate works. In our case, the fire safety requirements affect the control of system conditions, affecting all other parameters $x_i(t)$ (1). We check the completeness of the system description with the selected technique. If the syntax of the control method satisfies our requirements, then go to step 3, otherwise we return to the selection of a new control method step 2. If the method for control is not found, but the description and control of the system only by the name of the variables $x_i(t)$ satisfies us, then go to step 3, otherwise 1 step.

3 step. We check the control for optimality. If the solution is not optimal, according to the set objective functions, and does not suit us, then we return to step 1, otherwise we proceed to step 4. 4 step. We evaluate the effectiveness of strategy. If efficiency does not suit us, then we set new control actions and go to step 1, otherwise the end of the algorithm.

**4. Practical task**

An example 1. A woodworking enterprise [17]. Calculations of $G$ (4) were performed in the author's software package. Dimension $n$=1.2 million, $T$=5.5 years, $t$=1 month. In $t$=29 $J$ - an increase in forest processing from 1 to 2 million cubic meters per year ($W$). The incoherence of business processes among themselves, the presence of gaps in the volume of processing of semi-finished products, products in the material flows of processes, the lack of clear forms of document control procedures forced new investors (owners) to perform BP reengineering. The enterprise has completed a description of business processes (BP) and the implementation of a quality management system (TQM) in accordance with GOST R ISO 9001:2001, the structure of which is given by the matrix $B(t)$. In the framework of this decision system $u(t)$ is accepted by the owners (heads of departments, functional areas) of business processes in the form of planned values of business process performance indicators $x(t)$ taking into account the requirements of the accounting policy of the enterprise $A(t)$. Approval of the parameters of business processes remains with the head of the organization. The planned and actual parameters of the system are fixed through the integral indicator $G(t)$ according to formula. The figure shows the operating modes of the enterprise in a) standard mode and b) after the description of business processes and the implementation of the TQM. The analysis revealed that the processes were performed, but were not displayed in accounting and in the regulations. The total value of the non-visible material flows of the production processes of "shadow" business processes amounted to 4.5 billion rubles in year. A more detailed analysis of the data can be found in a separate work [18].

An example 2. The woodworking enterprise [17]. At the same enterprise, after the business process reengineering performed in the example 1 above, the PMBOK [16] process control standard was introduced into the matrix of the control structure $B(t)$. Within this system, control decisions $u(t)$ are made by the owners (heads of departments, functional areas) of business processes in the form of planned values of business process performance indicators $x(t)$ in the context of knowledge areas and gradation of control processes according to the PMBOK standard including the requirements of the accounting policies of the enterprise $A(t)$. The approval of all parameters of the business processes remains with the head of the organization. The planned and actual parameters of the system are fixed through the integral indicator $G(t)$ according to the formula: a figure 1 shows that the integral indicator $G$ (implementation of PMBOK) takes into account the improvement in the preparation processes, from 13 to 19 periods, by the time production increases from the 30 period. The figure 1 shows the implementation of the PMBOK standard itself is difficult to estimate in terms of value, but the integral indicator reflects the change in layers (process) in the system and their connectivity through the mode of operation c). Due to the introduction of PMBOK and the streamlining of controlled processes, unnecessary business processes that interfere with the increase in production were removed, therefore, there was a decrease in costs for this preparatory stage. A more detailed analysis of the data can be found in a separate work [19].

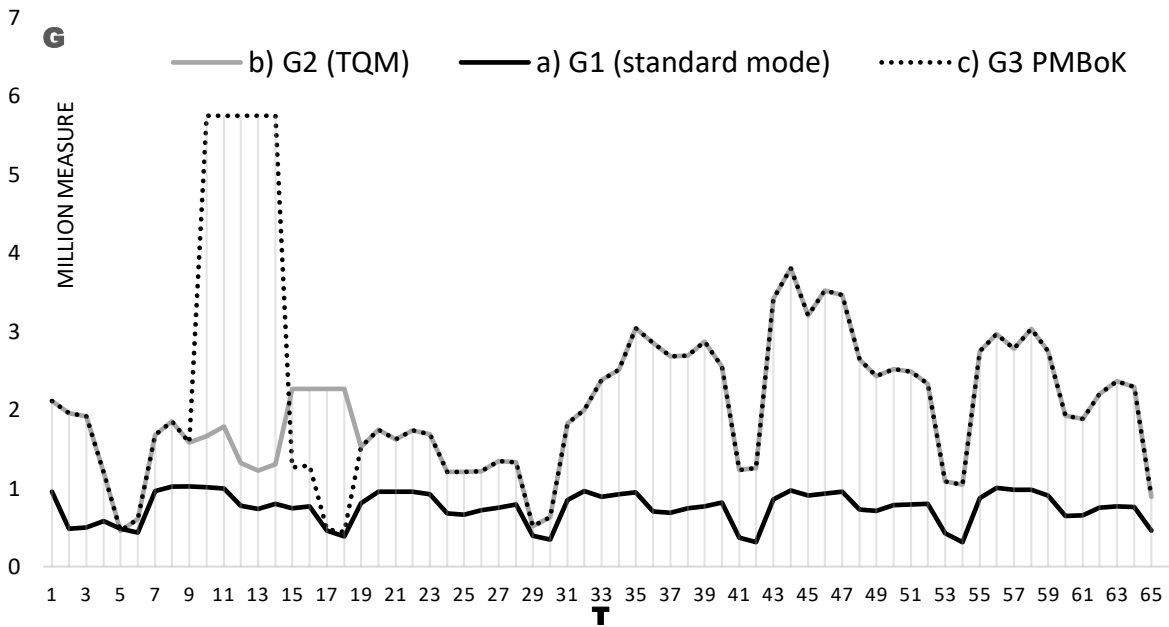

**Figure 1.** The values of the integral exponent *G* of the control regime a), b), c).

An example 3. A woodworking enterprise. A resident of the SEZ is a construction company. Dimension $n=400$, $T=5$ years, $t=1$ month, $J$ – two different staff strategies that differ in the set of job responsibilities involved, $W$ – the strategy of human resources management (HR). In the third example, a construction company established in 2002 is used as an economic object for research. The company carried out the construction of 4 residential buildings. It was monitored from December 2003 until mid-2008 (52 periods). The need to use this enterprise, as an example, is due to the fact that it is necessary to record changes in the state of the system from changes in the parameters characterizing the result of applying control actions to the functional responsibilities of personnel (HR) in a well-studied, real-life economic facility. Simply put, already accepted control actions are not reviewed so as not to mix the system response to them with the studied control actions (functional responsibilities of the staff). The figure 1 shows that in normal mode, the system has an initial set of fulfilled functional duties of personnel from job descriptions characterizing HR ($G_1$- the basic strategy), and the indicator value is 153,080. The restrictive regime (infection of personnel with the Covid-19 coronavirus) was introduced by excluding engineering personnel from activities important for construction. The list of functional duties of personnel blocked by Covid-19: design engineer in the period 1-19 - providing technical documentation for the construction site and in the periods 38-42, 50-54, 62-66 - development of technical specifications; engineer of the technical department in the period 1-20 - development of technical solutions, in the period 28-36 - development of projects for the production of works; design engineer in the period 1-27 - participation in the installation of structures, testing and commissioning, in the period 28-70, in the period 71-73 - participation in the installation of structures. A table 1 shows the entire list of functional (staff responsibilities) aren't performed by personnel due to the Covid-19 coronavirus (stop business-process).

**Table 1.** Incomlete job descriptions due to Covid-19 ($G_1$).

| Period (t) | Employee's position | Functional (staff responsibilities) $x_i$ |
|---|---|---|
| 1 - 19 | Concept engineer | Providing technical documentation for the construction site |
| 1 - 20 | Proof engineer | Creation of technical solutions |
| 28 - 36 | Proof engineer | Creation of work production projects |
| 38 – 42 | Proof engineer | Creation of work production projects |
| 50, 51 | Proof engineer | Creation of work production projects |
| 55 – 73 | Proof engineer | Creation of work production projects |
| 1 - 27 | Facility design engineer | Control of installation of structures, testing and commissioning of facilities |
| 28 - 70 | Facility design engineer | Designing the main sections of the project |
| 71 - 73 | Facility design engineer | Control and participation in the installation of structures |
| 1 - 32 | Comprehensive services | Implementation of construction supervision, participation in seminars and conferences |
| 50 – 51 | Comprehensive services | Author's supervision |
| 62 - 63 | Comprehensive services | Author's supervision |
| 1 – 18 | Project chief engineer | Project Coordination |
| 66 - 73 | Project chief engineer | Project Coordination |
| 38 - 39 | Project chief engineer | Author's supervision |
| 1 - 20 | Section foremaster | Implementation of technical control over the implementation of construction and installation work; acceptance of completed volumes |
| 25, 26 | Section foremaster | Material Accounting |
| 38 – 39 | Section foremaster | Material Accounting |
| 50 – 51 | Section foremaster | Material Accounting |
| 62 - 63 | Section foremaster | Material Accounting |
| 1 - 20, 25, 38, 50, 62 | Site supervisor | Performing construction and installation work |
| 1 - 20 | Production and Technical Department Engineer | Determination of the volume of work performed |
| 1 - 20 | Proof engineer | Keeping records of completed construction and assembly works and reporting on the implementation of plans |
| 1 - 20 | Occupational safety engineer | Monitoring the implementation of labor protection measures |

Blocked functional responsibilities of engineering personnel affect the performance $x_i$ of system functions. A figure 2 shows the estimate of HR in the infection regime ($G_3$ – strategy 3) is 155,150. Calculation result: the basic strategy ($G_1$=155,896), the strategy 3 ($G_3$=155,150). Functions blocked by Covid-19 can be restored within 1 months by attracting the services of third-party organizations in the amount of 188.8 million rubles, which is equivalent to a 31% rise in the financial cost of the entire project for 4 years.

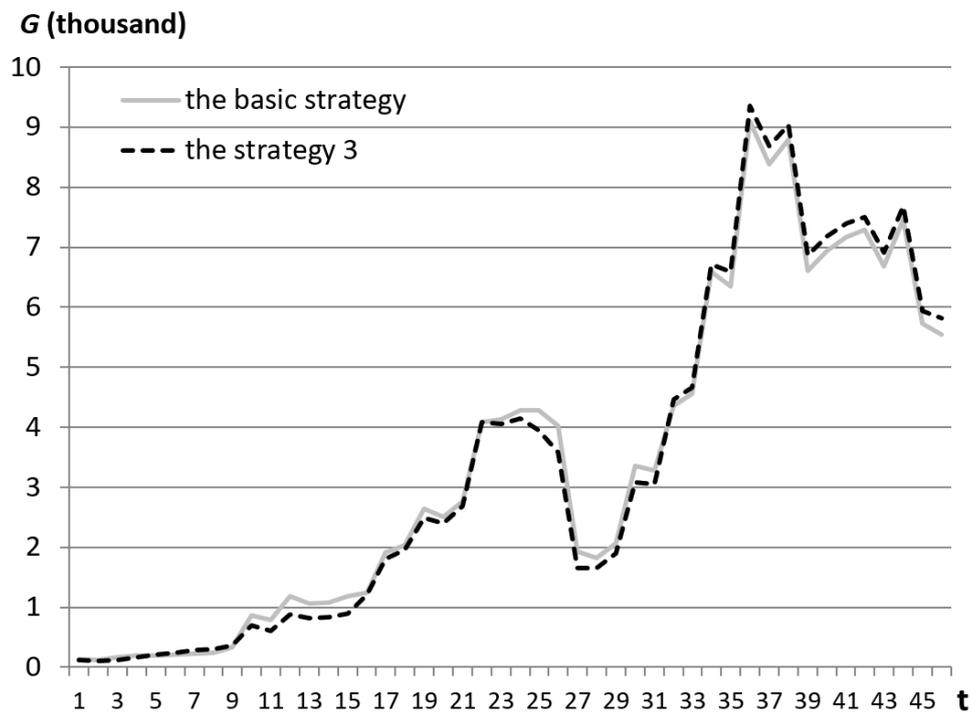

**Figure 2.** Parameter dynamics of $G$.

Then the influence of stop business-processes on the system is measured as the difference between the state of the control system parameters $\Delta G = G_3 - G_1$ and is equal to minus 745.9. The value is due to the higher multi-factorial response of the system to more difficult conditions for it. The multifactorial effect in systems is similarly disclosed in a separate research [5].

## 5. Conclusion
An assessment is made by an integral indicator of a dynamic system with a variable dimension at different time periods under the influence of internal factors and environmental parameters. The research showed significant changes in the values of the integral indicator characterizing the state of a dynamic system when implementing control standards. In the first and second examples, the underestimation of the hidden material flows of production processes of "shadow" business-processes turned out to be 4.5 billion rubles a year. In the third example, the loss of the enterprise from stop business-processes for 4 years will amount to 188.8 million rubles. The purpose of the research has been achieved.